\documentclass[11pt]{article}

\usepackage{amssymb,amsmath, epsfig}
\usepackage{changepage}
\usepackage{caption}
\usepackage[section]{placeins}
\makeatletter
\AtBeginDocument{%
  \expandafter\renewcommand\expandafter\subsection\expandafter{%
    \expandafter\@fb@secFB\subsection
  }%
}
\makeatother

\def\p{\partial}

\def\bbb{\noindent}  
\def\u{\vskip  .075 in}
\def\nh{\noindent\hangindent=1 true cm \hangafter = 1}

\def\EX{{\rm E}}
\def\nh{\noindent\hangindent=1 true cm \hangafter = 1}
\def\C{{\rm Cov}}
\def\u{\vskip  .1 in}

\def\B {\begin{eqnarray*}}

\newcommand{\bel}[1]{\begin{equation}\label{#1}}

\newcommand{\be}{\begin{equation}}
\newcommand{\qe}{\end{equation}}
\newcommand{\ee}{\end{equation}}
\newcommand{\baS}{\begin{eqnarray}}
\newcommand{\ba}{\begin{eqnarray}}
\newcommand{\ea}{\end{eqnarray}}

\def\EN{\end{eqnarray*}}
\begin{document}
\title{Numerical solutions of some hyperbolic stochastic partial differential equations with mixed derivatives including sine-Gordon equation\\
\   \\
{\normalsize
Henry C. Tuckwell$^{1,2\dagger, *}$\\   \
\  \\ 
 \              \\
$^1$ School of Electrical and Electronic Engineering, University of Adelaide,\\
Adelaide, South Australia 5005, Australia \\
 \              \\
$^2$ School of Mathematical Sciences, Monash University, Clayton,
Victoria 3800, Australia \\
 \              \\
%$^1$ Max Planck Institute for Mathematics in the Sciences\\
%Inselstr. 22, 04103 Leipzig, Germany\\
\     \\
\     \\
\vskip .25 in
%$^*$ Also at,   Max Planck Institute for Mathematics in the Sciences,
%Inselstr. 22, 04103 Leipzig, Germany\\
\  \\
$^{\dagger}$ {\it Email:} henry.tuckwell@adelaide.edu.au\\ }}

\maketitle

%
%
%\newpage 
\begin{abstract}  
We consider linear and nonlinear hyperbolic SPDEs with mixed derivatives with additive space-time Gaussian white noise of 
the form $Y_{xt}=F(Y) +  \sigma W_{xt}.$  
Such equations, which transform to linear and nonlinear wave equations, 
including Klein-Gordon,  Liouville's and the sine-Gordon equation, are related to what Zimmerman (1972) called a diffusion equation.  An explicit numerical scheme is employed in both deterministic
and stochastic examples.  The scheme is checked for accuracy against
known exact analytical solutions for deterministic equations.  In the stochastic case with $F=0$, solutions yield sample paths for the Brownian sheet whose statistics match well exact values.
Generally the boundary conditions are chosen to be initial values $Y(x,0)$ and boundary values $Y(0,t)$ on the quarter-plane or subsets thereof, which have been shown to lead to existence and uniqueness of solutions. For the linear case solutions are compared at various
grid sizes and wave-like solutions were found, with and without noise, for non-zero initial and boundary conditions. Surprisingly, wave-like structures seemed to emerge with zero initial and bondary conditions
and purely noise source terms with no signal. 
Equations considered with nonlinear $F$
included quadratic and cubic together with the sine-Gordon equation.
For the latter,  wave-like structures were apparent with $\sigma \le 0.25$ but they tended to be shattered at larger values of $\sigma$. 
Previous work on stochastic sine-Gordon equations is briefly reviewed.

 \end{abstract}
\centerline {{\it Keywords}: SPDE, mixed partials, space-time white noise, sine-Gordon}

\newpage
\rule{60mm}{1.5pt}
\tableofcontents

\rule{60mm}{1.5pt}
%
%Note new reference.
%Essential Psychopharmacology: Neuroscientific Basis and Practical Applications
% By S. M. Stahl. 
\newpage
  \section{Introduction}
There have been many recent articles involving the applications of stochastic partial 
differential equations. These include parabolic equations which have enjoyed widespread attention in neurobiological
applications
%(Albeverio et al., 2011, 2013;
(Boulakia et al., 2014; D\"orsek et al, 2103;
Faugeras and MacLaurin, 2014; Khoshnevisan and Kim, 2015;
Petterson et al., 2014; Stannat, 2013; Tuckwell, 2013a, 2013b) and to a lesser
extent hyperbolic equations (Hajek, 1982). 
 These applications have often employed 
two-parameter Wiener processes, or Brownian motion,  $\{W(x,t), x \in X, t \in T\}$,
with mean zero and covariance ${\rm Cov} [W(x,s),W(y,t) ]=\min(x,y) \min(s,t)$ where $X$ and $T$ are
sub-intervals of  $R$ or $R^+$ (or their formal derivatives, space-time white noise
$\{w(x,t)\}$). Such two-parameter processes 
 first appeared over 50 years ago  
in works by Kitagawa (1951), $\check{\rm C}$encov,  (1956 ) and Yeh (1960).

Zimmerman (1972) constructed a stochastic integral with respect to $W$ and 
obtained a solution to what that author called a diffusion
equation,
\ba \Delta Y(x,t)&=& Y(x,t)-Y(x,0)-Y(0,t) +Y(0,0)            \nonumber \\
&=& \int_0^t\int_0^x m(u,v,Y(u,v))dudv + \int_0^t\int_0^x  \sigma(u,v,Y(u,v))
dW(u,v) \ea
where $m$ and $\sigma$ are  Baire functions satisfying suitable
growth conditions. It was shown, inter alia,  that under the stated conditions,
solutions $Y$ had sample functions which were almost all 
continuous
 and were uniquely determined.
 
Yeh (1981) considered solutions of a stochastic differential equation
on $D=[0,\infty] \times [0, \infty]$ 
written as
\be dY(x,t)=m(x,t,Y)dxdt + \sigma(x,t,Y)dW(x,t)\ee
with a boundary condition 
\be Y(x,t)=Z(x,t), \hskip .25 in  {\rm for} (x,t) \in \partial D \ee
where $\partial D$ is the boundary of $D$ and $Z$ is a random
process with continuous sample functions and locally bounded second
moment on  $\partial D$.  It was proved that under suitable growth
conditions on $m$ and $\sigma$ that  a strong solution existed with
pathwise uniqueness. Interestingly, it was stated that  the system is non-Markov (Yeh, 1981).

It seems that the stochastic equations (1) and (2) have not been
employed as written in any biological, engineering or physical science modeling although we shall see that they transform to some well-known equations
of mathematical physics.

\section{A simple stochastic PDE with mixed derivatives}
The integral equation (1) or the differential equation (2) can be formally
recast as the stochastic partial differential equation 
\be  \frac{\p ^2Y}{\p x \p t} = m(x,t,Y) + \sigma (x,t,Y) \frac{\p ^2W}{\p x \p t}. \ee
This is a hyperbolic equation which is similar  to the  general
class of stochastic partial differential equations 
 considered by Hajek (1982)
 \be \frac{\p ^2Y}{\p x \p t} - a(Y) \frac{\p Y}{\p x} \frac{\p Y}{\p t} -b(Y)
- c(Y)w(x,t)=0, \ee
where $w$ is a space-time white noise with mean zero and covariance function
\be \C [w(x,s), w(y,t)]= \delta(x-y) \delta(s-t). \ee
In Hajek's equation the drift and diffusion coefficients depend only on $Y$. Hajek 
also proved existence and uniqueness of solutions and described a Stratonovich type
calculus. 

%The boundary condition for (5) was stated as $Y(x,t)=0$ if $xt=0$, and it
%was proved that solutions had a generalized Markov property. 

Identifying $t$ as a time parameter and $x$ as a spatial parameter, in the rest of this article we are concerned with temporally and spatially homogeneous
equations with the structure 
\be  \frac{\p ^2Y}{\p x \p t} = F(Y) + \sigma (Y) \frac{\p ^2W}{\p x\p t}, \ee
which is obtained from Hajek's equation by putting $a=0$ and substituting
$F$ and $\sigma$ for $b$ and $c$, respectively.

\section{Linear examples with space-time white noise} 
A simple equation with the structure of Eq. (7) is the
linear stochastic PDE with $F(Y)=\alpha Y$ and $\sigma(Y)=\sigma$,
where $\alpha$ and $\sigma$ are real constants. Thus
\be  \frac{\p ^2Y}{\p x \p t} =\alpha Y  + \sigma \frac{\p ^2W}{\p x \p t}, \ee
or, using subscript notation for partial differentiation,
\be Y_{xt} = \alpha Y + \sigma W_{xt}. \ee
Generally, in physical or natural sciences, it would be of interest
to find solutions on $ [0,x_f] \times [0, t_f]$ where $x_f \leq \infty$ and
$t_f \leq \infty$ with given values of $Y$, possibly random, on
the boundaries.
 
When $\sigma=0$, $Y$ is deterministic and an exact general solution may
be written as the linear combination
\be Y(x,t) = c_1 \exp(\alpha x)\exp(t) + c_2\exp(x)\exp(\alpha t)\ee
$c_1$ and $c_2$ are arbitrary constants. Note that an arbitrary constant cannot be added
to this general solution and the boundary conditions giving values of
$Y(0,t)$ and $Y(x,0)$ are, if this form of the solution is used,
being determined by the values of $c_1$ and $c_2$, are not arbitrary. 

\subsection{Numerical scheme with $\sigma=0$}
To solve Eq. (8) without noise numerically on  $ [0,x_f] \times [0, t_f]$, discretize $x$ and $t$ with steps of $\Delta x$ and $\Delta t$, so that $x_i=(i-1)\Delta x, i=1,2,\dots, n_x +1$
and  $t_j=(j-1)\Delta t, j=1,2,\dots, n_t+1$, so that $ \Delta x = x_f/n_x$ and
$\Delta t = t_f/n_t$.  Then approximating $Y(x,t)$ at the grid point $(x_i, t_j)$ by 
$Y_{i,j}$ and using the definition of second order partial derivative leads to
the scheme, 
\be  Y_{i,j} = Y_{i-1,j} + Y_{i,j-1} + \Delta x \Delta t F(Y_{i-1,j-1}) - Y_{i-1,j-1}. \ee
Supposing the initial values (IC) at $t=0$ are
\be  Y(x,0)=f(x), 0 \le x \le x_f < \infty,  \ee
and the boundary values (BC) at $x=0$ are 
\be Y(0,t)=g(t), 0 \le t \le t_f < \infty,  \ee
with the necessary consistency condition $f(0)=g(0)$. 
Then the values of $Y_{i,j}$ on the boundaries are
$Y_{i,1}=f(x_i), i=1,\dots, n_x+1$ and $Y_{1,j}=g(t_j), j=1,\dots, n_t+1$.
Such forms of initial and boundary values for linear and
nonlinear PDEs of the form
$Y_{xt}=F(Y)$ have been employed by Fokas (1997), Leon and Spire (2001), Leon (2003) and Pelloni (2005), mainly for the sine-Gordon equation. Pelloni (2005) pointed out that the use of such IC/BC
gives a well-posed problem and no additional constraints are 
required or necessary to define a solution. In fact, the
solution at $(x_0,t_0)$  depends only on solution values at
$x<x_0$ and $t<t_0$ so that one may continue to integrate
the PDE indefinitely away from the IC/BC conditions (12) and (13). 
Given  $Y_ {i,1}, \forall i$ and $Y_{1,2}$ enables the calculation of
$Y_{i,2}, \forall i$. Similarly, given $Y_{1,j} \forall j$ and $Y_{2,1}$ enables
the calculation of $Y_{2,j}, \forall j$. Proceeding in this fashion, $Y_{i,j}$ is
determined $\forall i,j$ from the initial/boundary values. 

     \begin{figure}[!h]
\begin{center}
\centerline\leavevmode\epsfig{file=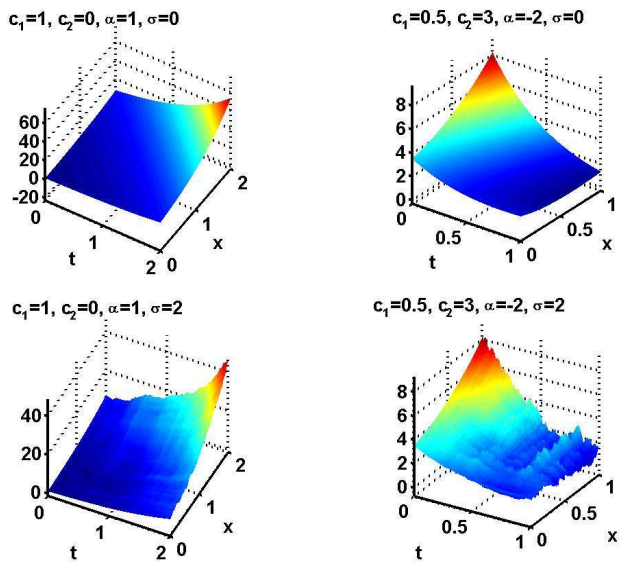,width=4.25in}
\end{center}
\caption{The top two panels show the numerical solution of the linear
PDE (8) without noise and with the parameters stated on the
figure.  The lower two panels show sample paths obtained by
numerical solution of (8) with the same parameters but including noise with $\sigma=0.5$ on the left and $\sigma=0.1$ on the right. 
 } 
\label{fig:wedge}
\end{figure}

Exact deterministic solutions, that is, for $\sigma=0$, were checked
against the solutions of the form of Equ. (10).  For a first example,
putting $c_1=1$, $c_2=0$ and $\alpha=1$ on  $ [0,2] \times [0, 2]$
with $n_x=n_t=100$ or $\Delta x=\Delta t = 0.02$  gave $Y(2,2)=53.290$ 
which compares with the exact value of $e^4=54.598$.
However, with  $n_x=n_t=500$ or $\Delta x=\Delta t = 0.004$,
the numerical value was 54.33, being much closer to the exact solution.
The numerical solution with  $n_x=n_t=100$ is shown as the surface
in the top left part of Figure 1. 
As a second example, take $c_1=0.5$, $c_2=3$ with $\alpha=-2$ on
 $ [0,1] \times [0, 1]$
with $n_x=n_t=100$ or $\Delta x=\Delta t = 0.01$. The numerical scheme
 gave $Y(1,1)=1.2728$ compared with the exact value $3.5e^{-1} \approx 1.2876$. 
With  $n_x=n_t=500$ or $\Delta x=\Delta t = 0.002$, the numerical scheme 
gave $Y(1,1)=1.2846$.  The maximum value of $Y$ over $ [0,1] \times [0, 1]$
was 8.2225 for both grid sizes. The numerical solution is shown for the
grid with $n_x=n_t=100$ in
the top right part of Figure 1. 

%     \begin{figure}[!h]
%\begin{center}
%\centerline\leavevmode\epsfig{file=FIM4.eps,width=5.5in}
%\end{center}
%\caption{Illustrating the accuracy of the numerical solution of (8)
%computed on  $ [0,2] \times [0, 3]$ .
%Here $\alpha=0.5$, $c_1=0.5$, $c_2=1$. The deterministic exact
%solution is shown in the top panel and the relative error
%of the numerically computed solution for a grid 
%with $n_x=1000$, $n_t=1200$ is shown in the lower panel.
%The maximum relative error is 0.0016.  
% } 
%\label{fig:wedge}
%\end{figure}

            As a further test of the numerical scheme, the  same equation
was considered on  $ [0,2] \times [0, 3]$ with parameters $\alpha=0.5$, $c_1=0.5$, $c_2=1$. For a grid with $n_x=n_t=500$, the maximum
relative error was 0.0035, and with $n_x=700$, $n_t=800$ it was
0.0023. For a still finer grid  $n_x=1000$, $n_t=1200$ the maximum relative error was 0.0016. For this latter case the exact solution 
\be Y(x,t)= 0.5\exp(0.5x)\exp(t) + \exp(x)\exp(0.5t) \ee
evaluated at grid points and called $Z(i,j)$ is shown in the top part of Figure 2, the bottom part giving the relative error against the exact solution. As the mesh becomes finer, the relative error diminishes at all
grid points. The maximum value of 0.16\% for the case shown
points to the accuracy of the numerical scheme. 
\subsection{Stochastic case}
To simulate the additional white noise term in Equation (9),
we approximate  $\sigma \frac{\p ^2W}{\p x \p t} $ as in the numerical solution
of the spatial Fitzhugh-Nagumo equation (Tuckwell, 2008) with space-time white noise, where results were checked against analytical
solutions for the moments.  Thus at $(x_i,t_j)$   
\be \frac{\p ^2W}{\p x \p t}  \approx \sigma \frac{1}{\sqrt{\Delta x \Delta t}}
N_{ij}  \ee
where 
where the $N_{i,j}$ are independent standard (zero mean, unit variance)
normal random variables, which will be generated by computer random number
generator. Thus the discretized version of (8) is
\be  Y_{i,j} = Y_{i-1,j} + Y_{i,j-1} + \Delta x \Delta t F(Y_{i-1,j-1}) - Y_{i-1,j-1} + \sigma \sqrt{\Delta x \Delta t} N_{ij}. \ee 

%the following term is added to the right hand side of Equation (11)
%\be \sigma w(x,t) \approx \sigma \sqrt{ \frac{\Delta t}{\Delta x}} N_{ij} \ee
%where $N_{ij}, i=1, \dots, n_x+1, j=1, \dots, n_t+1$ are independent
%standard normal random variables.
A sample path with the inclusion of noise with $\sigma=2$ for the first
deterministic example ($c_1=1, c_2=0, \alpha=1$) is shown in the bottom left panel of Figure 1 (below the deterministic solution).  Here again $n_x=n_t=100$ and $x_f=t_f=2$.
For 100 trials, the standard deviation of 
$Y(2,2)$ was 12.75 and the mean of $Y(2,2)$ was 54.66, compared to the no noise numerical value of 53.29 with
the same grid size  and the exact value of 54.60. 

In the bottom right panel of Figure 1 is shown a sample path,
 with $\sigma=2$, for the second deterministic example
($c_1=0.5, c_2=3, \alpha=-2$), again on  $ [0,1] \times [0, 1]$ with
$n_x=n_t=100$. For 2 sets of 100 trials, the standard deviation of $Y(1,1)$
was 1.406 and 1.265, with means of 0.976 and 1.20, respectively,
the latter results being considerably less than the exact mean
of 1.2876. When the grid was made finer with $n_x=n_t=200$,
the mean was closer to the exact result at 1.26 and the standard deviation was 1.477 (more details below). 
\u
\subsection {Brownian sheet sample paths}
\u
\bbb It was of interest to use the numerical scheme of Eq. (16)  
with $F=0$ together with boundary values of 0 to generate
sample paths for the two-parameter Wiener process or Brownian
sheet (Welner, 1975; Adler, 1978; Koshnevisan, 2001). This was done
with $\sigma=3$ on $ [0,1] \times [0, 1]$ and on $ [0,2] \times [0, 2]$
for two different mesh sizes, $n_x=n_t=100$, $n_x=n_t=400$,  
with 500 and 200 trials respectively. Examples of the sample paths
are shown in Figure 2, the top two panels being for $W$ on
 $ [0,1] \times [0, 1]$ 
and the bottom two on $ [0,2] \times [0, 2]$.

Table 1 gives the values of the mean  $\EX [W(x_f,t_f)]$ and the standard deviation $SD [W(x_f,t_f)]$ 
of $W(x_f, t_f)$ from the simulations. The value of the exact
mean is $\EX [W(x_f,t_f)]=0$ and that of the standard deviation
is $SD [W(x_f,t_f)]= x_f \sigma$, being either 3 or 6.  The approximate 
95\%
confidence intervals for the mean are given in column 5.
For both sets of trials on $ [0,1] \times [0, 1]$  the sample
mean is well within the 95\% confidence limits, but in the first
run with 200 trials with the finer mesh on $ [0,2] \times [0, 2]$,
the sample mean is just outside the 95\% confidence interval. In a second
set of 200 trials, the sample mean was well within the confidence limits.

\begin{center}

\begin{table}[h]
    \caption{Numerical results for simulation of Brownian sheet: in all cases $\sigma=3$}
\smallskip
\begin{center}
\begin{tabular}{cccccc}
\hline 
$n_x=n_t$   & No. Trials  & $x_f=t_f$ & $\EX[W(x_f,t_f)]$  & 95\% C.I. for \EX & $SD[W(x_f,2´t_f)]$ \\
  \hline
100 & 500 & 1 & -0.165 & $\pm 0.263$  &  3.104 \\
400 & 200 & 1 & 0.237 &  $\pm 0.416$ & 3.296 \\
100 & 500 & 2 & 0.041 & $\pm 0.526$ & 5.803 \\
400 & 200 & 2 & -0.893 & $\pm 0.832$ & 5.867 \\
400 & 200 & 2 & 0.310 & $\pm 0.832$ & 6.589 \\

      \hline
\end{tabular}

\end{center}

\end{table}

\end{center}

      \begin{figure}[!h]
\begin{center}
\centerline\leavevmode\epsfig{file=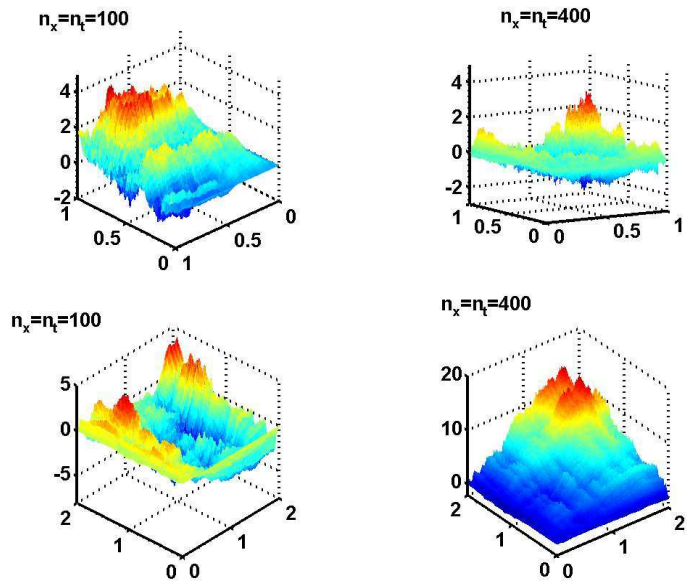,width=4.5in}
\end{center}
\caption{Four simulated sample paths of a Brownian sheet $W$
generated by applying Eq. (16) with $F=0$ and $\sigma=3$. 
The top two panels are on $ [0,1] \times [0, 1]$  whereas the bottom two
are 
on $[0,2] \times [0,2]$. The grid is finer (400 by 400) in the
left hand cases, being 100 by 100 in the right hand cases.  
 Boundary values $W=0$ are applied on the $x$ and $t$ axes. }
%on $\{x \in [0,X] , t=0\}}$ and $\{x=0, t \in [0,T]\}$}.  
\label{fig:wedge}
\end{figure}

\subsection{Further numerical results for the linear stochastic Equation (8)}

Considering further the SPDE (8) with noise, it is of 
interest to see how the results might depend on grid size and 
number of trials.
To this end, with $c_1=1, c_2=0$ and $\alpha=1$,  several 
runs were done on  $ [0,2] \times [0, 2]$  as described in Table 2 with
the stated results for the mean and standard deviation
of $Y(2,2)$, being the maximum point in the deterministic case.
In general it is apparent that the value of $Y(2,2)$ is underestimated
in all the deterministic cases, though the numerical value seems
to approach the exact value as the grid becomes finer. Thus the
error with a grid of 100 x 100 gives an error of 2.39 \% whereas
a grid of 1000 x 1000 gives an error of 0.25 \%. 

 As it is apparent (and see also Section 3.5), 
by taking expectations and integrating directly in Equ. (8), that the mean of $Y$ with the additive white noise is the same as the deterministic
value, then the same remarks apply to the mean in the stochastic case,
as the difference between the exact
deterministic value and the mean with noise at $\sigma=0.1$ 
 drops (monotonically for these simulations) from 2.36 \% for a grid of 100 x 100 to
0.63 \% for a grid of 400 x 400. From the results of the last 4 rows
in Table 2 where $\sigma=0.5$, there is no consistent pattern 
of changes in the mean relative to the deterministic value as the
number of trials increases, though the set of results is very small.

\begin{center}

\begin{table}[h]
    \caption{Results for Equation (8) with $\alpha=1$ for various grid sizes and noise levels}
\smallskip
\begin{center}
\begin{tabular}{ccccc}
  \hline
    Run     & $\sigma$  & $n_x=n_t$ & \EX[Y(2,2)]  & SD[Y(2,2)] \\
  \hline
Exact & 0 &  - &   54.598   & 0 \\
\hline
Numerical deterministic & & & & \\
 & 0 & 100 & 53.290 & 0 \\
& 0& 200 & 53.936 & 0 \\ 
 & 0 & 500 & 54.331 & 0\\
 & 0 & 1000 & 54.464 & 0 \\
\hline
Numerical SPDE & & & & \\
No. trials & & & & \\
200 & 0.1 & 100 & 53.310 & 0.601 \\ 
200& 0.1 & 200 & 53.910 & 0.640 \\
200 & 0.1 & 300 & 54.133& 0.657 \\
200 & 0.1 & 400 & 54.253 & 0.645\\
100 & 0.5 & 200 & 53.787 & 3.112 \\
200 & 0.5 & 200 & 53.520 & 3.172 \\
300 & 0.5 & 200 & 54.166 & 3.129 \\
400 & 0.5 & 200 & 53.918 & 3.154 \\
  \hline
\end{tabular}

\end{center}

\end{table}

\end{center}

Further results were obtained for the stochastic PDE (8), 
  integrated with various noise levels $\sigma$ from 0 to 3, using the same parameters  as above but that in every case $n_x=n_t=200$ and
200 trials. When $\sigma=2$, for example, the mean and standard deviation of $Y(2,:)$  were plotted against time for $t$ up to 2. The agreement between the mean and the exact
deterministic value  $Y(2,:)=\exp(2)\exp(t)$  was excellent. The standard
deviation grew in an approximately linear fashion over the same time interval.  When the standard deviation of the solution value $Y(2,2)$ is plotted against
noise level, with all other parameters fixed, the growth of the standard 
deviation $ SD[Y(2,2)] $ is about linear with $\sigma$ and in fact a good fit to the
curve is given by 
\be  SD[Y(2,2)]  \approx 6\sigma, \ee
which can be compared with the value $2\sigma$ for the Brownian sheet. 
%      \begin{figure}[!h]
%%\begin{center}
%%\centerline\leavevmode\epsfig{file=FIM3.eps,width=5.5in}
%%\end{center}
%\caption{In the upper part the mean (red curve) and standard deviation 
%(blue curve) of $Y(x,t)$
%at $x=2$ as functions of time are shown
%from the numerical solution of $Y_{xt}=\alpha Y + \sigma W_{xt}$
%with $\alpha=1$, $\sigma=3$. Parameters are $c_1=1, c_2=0$
%$n_x=n_t=200$ on  $ [0,2] \times [0, 2]$ : 200 trials. Also shown is the
%exact deterministic  ($\sigma=0$) solution $Y(2,:)=\exp(2)\exp(t)$ (green curve).
%In the lower part the standard deviation of $Y(2,2)$ is plotted against
%the noise amplitude $\sigma$. Parameters otherwise the same as
%in the top part of the figure.}
%%on $\{x \in [0,X] , t=0\}}$ and $\{x=0, t \in [0,T]\}$}.  
%\label{fig:wedge}
%\end{figure}

\subsection{Other boundary conditions: waves and 
purely noise generated waves}
In the remainder of this section we further consider mumerical
integration of the 
linear equation 
 \be Y_{xt}=\alpha Y + \beta + \sigma W_{xt}. \ee
Because the noise term contributes zero
to the mean we expect that $\hat{Y} = E[Y]$ will satisfy 
the same equation as the deterministic solution, so that
\be \hat{Y}_{xt} = \alpha \hat{Y}  + \beta. \ee
With the parameters $\alpha=-1$ and $\beta=0$, the veracity
of this claim was tested on the square $ [0,5] \times [0, 5]$
with the initial/boundary conditions
\be Y(x,0)=1, 0 \leq x \leq 5, Y(0,t)=1, 0 \leq t \leq 5, \ee 
 and with a grid of 500 x 500. 
In Figure 3 are shown results for the deterministic solution
(red surface) and the means for 50 trials with $\sigma=0.05$ (blue surface) and with $\sigma=0.1$ (green surface). 
It can be seen that the three surfaces are almost identical
so that the mean can, with suitable geometry, mesh sizes and parameter sets
be found as the deterministic solution. This also provides
a heuristic test of the accuracy of a numerical scheme.

      \begin{figure}[!h]
\begin{center}
\centerline\leavevmode\epsfig{file=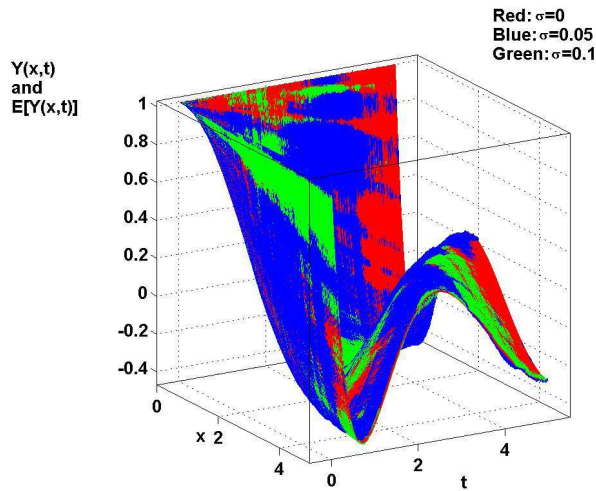,width=4.0in}
\end{center}
\caption{For $Y_{xt}= -Y + \sigma W_{xt}$ the deterministic
solution is shown in red and the means for the stochastic cases
$\sigma=0.05$ and $\sigma=0.1$ are shown in blue and green
respectively. The three surfaces are practically the same.
 Initial/boundary values both unity. }
%on $\{x \in [0,X] , t=0\}}$ and $\{x=0, t \in [0,T]\}$}.  
\label{fig:wedge}
\end{figure}

\subsubsection{Previous boundary condition}
With boundary conditions of the form of (10), wave-like solutions
were not apparent with the parameters $\alpha=-1$, $c_1=0.5$,
$c_2=1$, $\sigma=0$.  This was the case no matter how large a 
space- time
interval was considered, up to $ [0,20] \times [0, 20]$
with a grid $500 \times 500$.  Furthermore, the numerical solution
agreed precisely with the analytical solution.

\begin{figure}[!ht]
\begin{center}
\centerline\leavevmode\epsfig{file=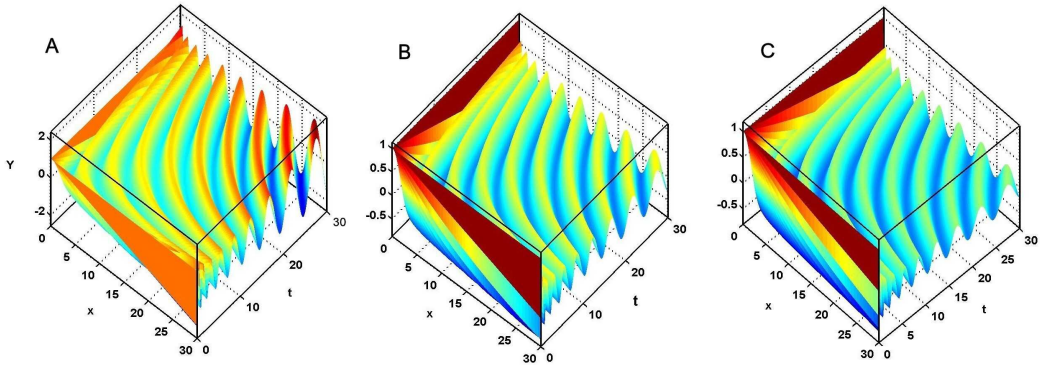,width=5.5in}
\end{center}
\caption{Numerical solution of  Eq. (18)  for various grid sizes  without
noise, $\sigma=0$,  and with IC/BC 
$Y(x,0)=1, 0 \leq x \leq 30$ and $Y(0,t)=1, 0 \leq t \leq 30$.
Parameters $\alpha=-1$, $\beta=0$. A,  grid $300 \times 300$,
B, $700 \times 700$, C, $1100 \times 1100$.} 
\label{fig:wedge}
\end{figure}
  
\subsubsection{Waves with non-zero boundary conditions}
With the same parameters as in the previous example
but with the initial/boundary conditions 
\be Y(x,0)=1, 0 \leq x \leq 30, Y(0,t)=1, 0 \leq t \leq 30, \ee 
with no noise ($\sigma=0$), waves of a sinusoidal shape
 form as depicted in the three 
panels of Figure 4.  
Solutions were obtained for meshes with $n_x=n_t= 300$, (left),
$n_x=n_t= 700$, (middle),  and  $n_x=n_t= 1100$, (right). 
The amplitude of the waves diminishes as the grid becomes finer,
presumably reflecting greater accuracy of the solution. 

This was examined further by computing the solutions at three different
mesh sizes ($500 \times 500$, $1000 \times 1000$ and $2000 \times 2000$) on a larger area,  $ [0,40] \times [0, 40]$. Results are shown in Figure 5, where $Y(x,40)$ is plotted against $x$. In each case there are
12 peaks but as the grid becomes finer, and presumably the solution
more accurate, the amplitude of the waves gets smaller and seems
to attain a constant value, along with a fairly constant spatial and temporal frequency.

The  appearance of wave-like solutions is not surprising since under the
coordinate transformation
\be \xi=x+t   \ee
\be \eta = x-t   \ee
the equation 
\be Y_{xt}=F(Y)  \ee
becomes the general wave equation
\be Y_{\xi \xi} - Y_{\eta \eta} = -F(Y) \ee which, in the case of a nonlinear
$F$ is called a nonlinear Klein-Gordon equation or with linear $F$ simply
the Klein-Gordon equation of mathematical physics. The linear case
is also a small amplitude approximation to the sine-Gordon equation
considered below. 
\begin{figure}[!t]
\begin{center}
\centerline\leavevmode\epsfig{file=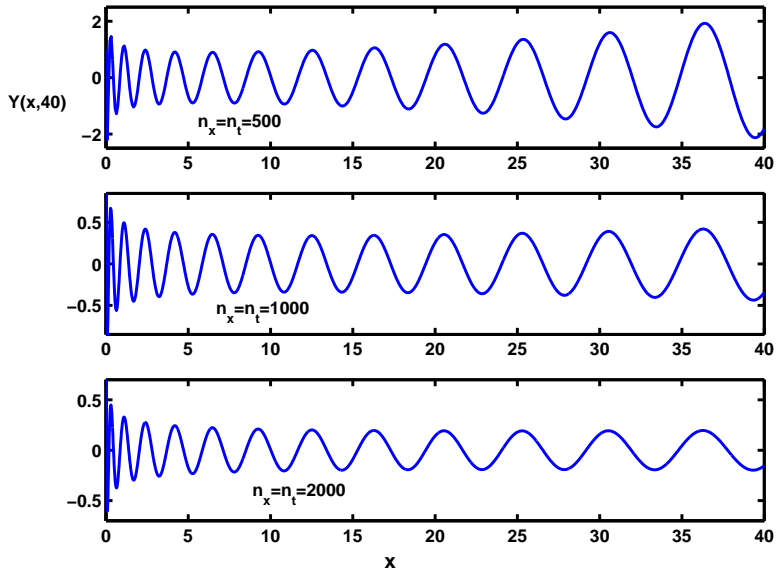,width=4.0in}
\end{center}
\caption{$Y(x,t_f)$ versus $x$  with $x_f=40$ and $t_f=40$ for the
solution without noise shown in Figure 5, but with meshes of various 
sizes. Top picture, $n_x=n_t=500$; middle, $n_x=n_t=1000$, 
and bottom $n_x=n_t=2000$.} 
\label{fig:wedge}
\end{figure}

When 2-parameter white noise with amplitude $\sigma=0.5$ was present, again with
$\alpha= -1, \beta=0$,  the mean and variance were obtained numerically
for 50 trials on  $ [0,30] \times [0, 30]$. The result for the mean $E[Y(x,t)]$  is plotted
in Figure 6, where again the grid is finer from left to right (left $300 \times 300$, middle $500 \times 500$, right $700 \times 700$). 
In the mean the waves are still discernible and the amplitude again
decreases as the grid becomes finer. There are 9 wavefronts discernible
in the right hand part of Figure 6 (stochastic case) and in the middle part
of Figure 4 (deterministic case), for which the grids are both 
$700 \times 700$. When these two plots are
rotated in order to make a comparison clearer,  for the stochastic case, the wave amplitude of $E[Y(x,t]$
increases as the wave progresses, but it is not known whether this effect
is real or due to an accumulated numerical error for large $x$ and $t$.

\begin{figure}[!ht]
\begin{center}
\centerline\leavevmode\epsfig{file=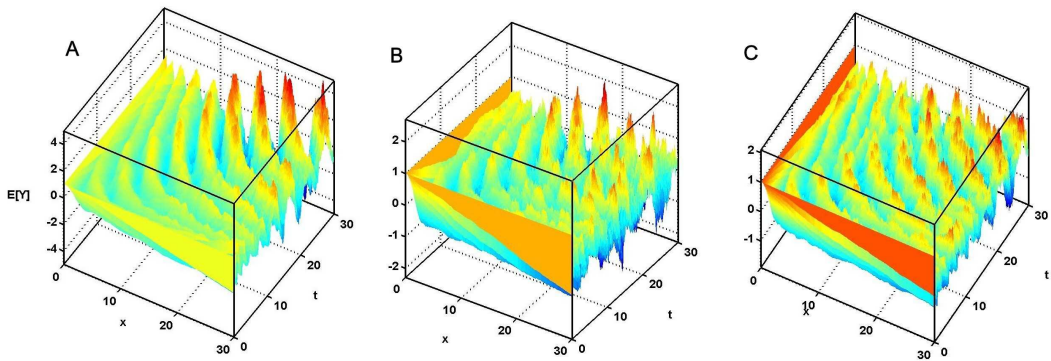,width=5.5in}
\end{center}
\caption{Mean of 50 trials for the solution of  Eq .(18) with
$Y(x,0)=1, 0 \leq x \leq 30$ and $Y(0,t)=1, 0 \leq t \leq 30$ with noise
 of amplitude $\sigma=0.5$. 
Parameters $\alpha=-1$, $\beta=0$. Grids:
A,  $n_x=n_t=300$, B,  $n_x=n_t=500$, 
C, $n_x=n_t=700$.} 
\label{fig:wedge}
\end{figure}

%
%\begin{figure}[!ht]
%\begin{center}
%\centerline\leavevmode\epsfig{file=F5CANGLE.eps,width=3in}\epsfig{file=F7AANGLE.eps,width=3in}
%\end{center}
%%\begin{center}
%%\centerline\leavevmode\epsfig{file=F7C.eps,width=3in}
%%\end{center}
%%\begin{center}
%%\centerline\leavevmode \epsfig{file=F7A.eps,width=3in}
%%\end{center}
%\caption{Highlighting the difference in amplitudes of waves in the
%deterministic case (left) and the mean for the stochastic case (right). These are 
%rotated 
%plots of the middle part of Figure 5 and the bottom part of Figure 7.} 
%\label{fig:wedge}
%\end{figure}

%\begin{figure}[!ht]
%\begin{center}
%\centerline\leavevmode\epsfig{file=FZWAVE1.eps,width=3in}\epsfig{file=FZWAVE1A.eps,width=3in}
%\end{center}
%%\begin{center}
%%\centerline\leavevmode\epsfig{file=F7C.eps,width=3in}
%%\end{center}
%%\begin{center}
%%\centerline\leavevmode \epsfig{file=F7A.eps,width=3in}
%%\end{center}
%\caption{Highlighting the difference in amplitudes of waves in the
%deterministic case (left) and the mean for the stochastic case (right). These are 
%rotated 
%plots of the middle part of Figure 5 and the bottom part of Figure 7.} 
%\label{fig:wedge}
%\end{figure}

\subsubsection{Wave-like behavior with noise and zero boundary conditions}
Having seen that solutions of $Y_{xt}=\alpha Y + \sigma W_{xt}$, with
$\alpha=-1$, $\sigma=0.5$ and with initial/boundary condition as in
(20)  or (21) gave wave-like solutions, it was decided to examine solutions
with the {\it zero} boundary conditions,
\be Y(x,0)= 0, \leq x \leq 20, Y(0,t)=0, 0 \leq t \leq 20.\ee 
With no noise the solution is not surprisingly identically zero
for all $x$ and $t$. 

With noise level set at $\sigma=0.5$ patterns emerge which suggest wavefront activity. The mean of solutions for 50 trials is shown in
Figure 7 for a grid 900 x 900. The pattern seems to self-organize into quite large wave-like segments with a wavelength which is fairly constant.

To the right of the  plot of $E[Y(x,t)]$ (Figure 7B) there is a plan view (looking in the
direction of -z) of $I_{[a, \infty]}(E[Y(x,t)])$ where $ 0 <a < \max (E[Y(x,t)])$ and  $I_A(z)$ is the indicator
function taking the value 1 if $z \in A$ and 0 otherwise.  That is, only the top portion of the surface lying above $z=a$ is shown. Such a plot
emphasizes the wave-like nature of $E[Y(x,t)]$. The value 
 of $a=0.15$. (In the caption of Figure 7 is also given an upper cut-off
value, but this is  above the maximum which makes it effectively infinite.) To further highlight the wave-like nature of 
(the mean of) the solutions, in Figure 7C a sequence of approximate
wavefronts is marked in the case of  a grid of
500 x 500.  The evidence for a wave-like structure seems quite
convincing. 

The appearance of these putative wave-like structures is unexpected.
There is no    ``signal''  because $\beta=0$ and the initial/boundary
conditions for $Y(x,t)$ are identically zero.  Hence the source
of the wave-like structures is purely space-time white noise acted upon
by the dynamical system. This response can be compared to the 
Brownian sheet sample paths with the same initial/boundary conditions
as depicted in Figure 2, where no distinct pattern of wavefronts is
apparent.  
It seems that the operator $\mathcal { L} $
defined through casting the SPDE as
\be \mathcal { L} \{Y\}  = \frac {\p^2 Y}{\p x \p t} + \alpha Y = \sigma W_{xt} \ee 
induces wave-like solutions from scattered local changes in $Y$
coming from the Brownian sheet.  The phenomenon which has
the appearance of producing organized clusters is reminiscent,
with some imagination,  of
theories of the creation of the  ``universe from nothing'' (Vilenkin, 1982, 1984; Krauss, 2012).

\begin{figure}[!ht]

\begin{center}
\centerline\leavevmode\epsfig{file=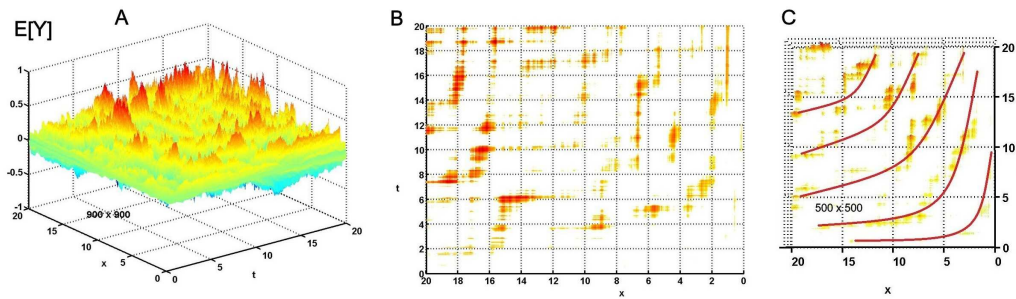,width=5.5 in}
\end{center}
\caption{Wave like structures with zero boundary conditions. 
A.  Plot of expected value of $Y(x,t)$ versus $x$ and $t$ for the  solution of  Eq. (18) with
$Y(x,0)=0, 0 \leq x \leq 20$ and $Y(0,t)=0, 0 \leq t \leq 20$ with noise
 of amplitude $\sigma=0.5$. 
Parameters $\alpha=-1$, $\beta=0$. Grid:  $n_x=n_t=900$. B.  On the Plan view with $E[Y(x,t)]$ cut off below a certain positive value to emphasize the wave-like structures,  $0.15 < E[Y] < 0.706=\max(E[Y])$.
C. Plan view for a coarser grid, 500 x 500,  with inserted 
curves which have been visually estimated to fit approximately to hypothesized wavefronts} 
\label{fig:wedge}
\end{figure}

%     \begin{figure}[!h]
%\begin{center}
%\centerline\leavevmode\epsfig{file=FZWAVE1APAT.eps,width=4.5in}
%\end{center}
%\caption{The same as the top right plot of the previous figure showing
%wave-like response of the system with zero signal and zero
%boundary conditions. Curves have been visually estimated to fit approximately to hypothesized wavefronts } 
%\label{fig:wedge}
%\end{figure}

\section{Simple nonlinear equations}
The above SPDEs have all been linear.  In this section we briefly consider
two nonlinear source functions $F$ which
are common in biological and physical models.

%In all cases the
%initial/boundary conditions are either (a) unity and (b) 0. 

\subsection{Quadratic $F$}
The classical example of a PDE with a quadratic source function is the parabolic equation
\be   Y_t = Y_{xx} + kY(1-Y) \ee
which was studied in the setting of spatial patterns of gene frequency by Fisher (1937) and Kolmogorov et al. (1937). 
In the setting of the present article we consider the SPDE
\be Y_{xt}=kY(1-Y) + \sigma W_{xt} \ee
which transforms to a nonlinear Klein-Gordon equation.

     \begin{figure}[!h]
\begin{center}
\centerline\leavevmode\epsfig{file=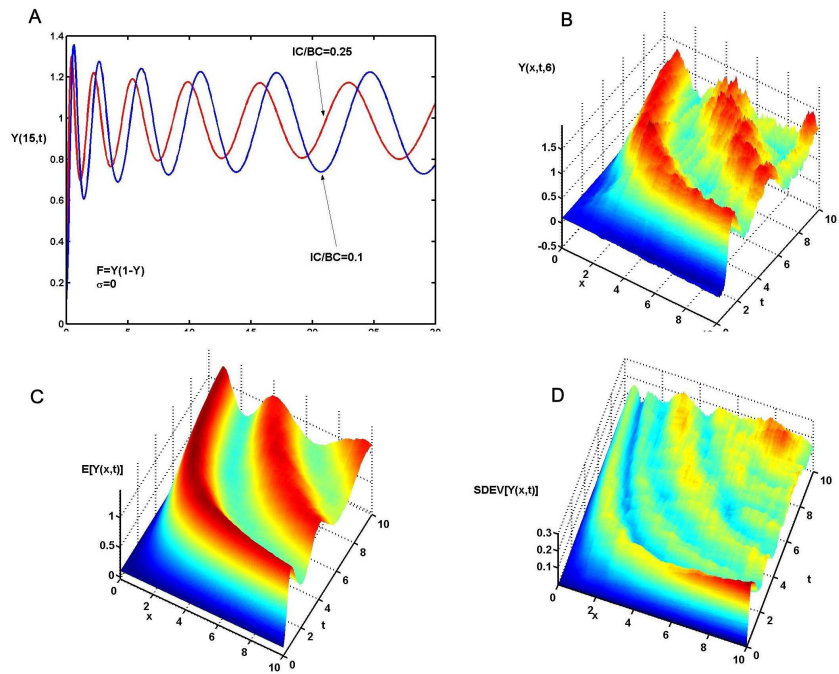,width=4.5in}
\end{center}
\caption{A. Values as functions of time of the numerical solutions of $Y_{xt}=Y(1-Y)$ on
 $ [0,30] \times [0, 30]$ at $x=15$. The blue curve is for the initial/boundary value of 0.1, whereas the red curve is for the initial/boundary value of 0.25.  B. A sample path for the SPDE (29) with $\sigma=0.06$  and IC/BC value 0.1.  This may be compared with the mean shown in the next diagram.  C and D. Mean (left) and standard deviation (right) 
obtained by  numerical solution (50 trials)
 of $Y_{xt}=Y(1-Y) +\sigma W_{xt}$ with $\sigma = 0.06$
and initial/boundary value of 0.1.} 
\label{fig:wedge}
\end{figure}
\subsubsection{Deterministic example}
For the deterministic equation with $\sigma=0$ numerical integration
was performed with the initial/boundary conditions of the form of (26)
but with $x_f=t_f=30$ and with the two values 0.1 and 0.25 rather than 0.
The results consisted of waves, similar to those shown in Figure 6,  which settled into an approximately
sinusoidal form. Sections are illustrated in Figure 8A where are shown plots of $Y$ at the middle $x$-value versus $t$, that is, $Y(t,15)$, the blue and red curves being for the initial/boundary values of 0.1 and 0.25
respectively.  The periods for both cases are similar but the 
amplitude is smaller for the larger initial/boundary condition.

%
%     \begin{figure}[!h]
%\begin{center}
%\centerline\leavevmode\epsfig{file=FQPATH1.eps,width=3in}\epsfig{file=FQPATH2.eps,width=3in}
%\end{center}
%\caption{Two sample paths for the SPDE (29) with $\sigma=0.06$ 
%and IC/BC value 0.1.  These may be compared with the mean shown in the next figure. 
%} 
%\label{fig:wedge}
%\end{figure}

\subsubsection{Stochastic example}
Numerical integration of the equation (29) was performed with many values of the parameters $x_f=t_f$,$n_x=n_t$ and $\sigma$, in conjunction with initial/boundary conditions 
\be Y(x,0)=0.1, 0 \leq x \leq x_f, Y(0,t)=0.1, 0 \leq t \leq t_f. \ee 
Well behaved solutions were only found for rather small values of
$\sigma$. 

The behavior of solutions on $[0,10]  \times [0,10]$
with a grid $n_x=n_t=500$ was investigated for various values
of $\sigma$. 50 trials were performed for each parameter set. 

A representative sample path (surface) 
is shown in Figure 8B with $\sigma=0.06$. The corresponding mean and standard deviation are given 
in Figures 8C and 8D, respectively.  The sample path does not depart greatly from the mean
but in others the second wave front has fractured (not shown). The wave-like behavior of the standard deviation roughly minics that of the mean, but there seem
to be additional weak wavefronts.

The most interesting aspect of solutions, however, is the extreme sensitivity to the magnitude of the noise parameter $\sigma$ near
a critical value. Table 3 gives maxima and minima for $E[Y]$ over the
region of integration starting at $\sigma=0.05$ and increasing
to $\sigma=0.06467150$. At the smallest of these values, the mean 
remains bounded for the whole region with a minimum value of 
0.1000 and a maximum of 1.3364. When $\sigma$ exceeds a 
critical value of (about) 
$\sigma_c=0.06467147$ the minimum of  $E[Y]$ becomes infinitely negative, 
but when $\sigma=\sigma_c$ the minimum of $E[Y]$ is finite and
positive at 0.0997.

\begin{center}

\begin{table}[h]
    \caption{Some solution properties for $Y_{xt}=Y(1-Y) + \sigma W_{xt} $}
\smallskip
\begin{center}
\begin{tabular}{lccl}
\hline 
$\sigma$   & $\max$ E[Y(x,t)] & $\min$ E[Y(x,t)] & Path in Figure 15 \\ 
  \hline
0.05 &  1.3364 & 0.1000  & blue dot-dash  \\
0.06 & 1.3131 & 0.0981 & blue solid \\
0.063 & 1.3078 & 0.0992 & red solid \\
0.064 & 1.3017 & 0.0995 & black solid \\
0.0645 & 1.3200 & 0.0993 & green solid \\
0.0646 & 1.3025 & 0.0989 & mauve solid \\
0.06465 & 1.3021 & 0.0996 & cyan solid\\
0.06467 & 1,3196 & 0.0996 & yellow solid\\
0.064671 & 1.3201 & 0.1000 & blue dashed \\
0.0646713 & 1.3109 & 0.0997 & red dashed \\
0.0646714 & 1.3258 & 0.995 & black dashed \\
0.06467145 & 1.3215 & 0.0995 & - \\
0.06467147 & 1.3140 & 0.0997 & - \\
0.06467148 & 1.3318 & $-\infty$ & - \\
0.06467150 &1.3095 & $-\infty$ & - \\
      \hline
\end{tabular}

\end{center}

\end{table}

\end{center}

The results for $E[Y(5,t)]$ for various  values of $\sigma$ from 0.05 to $\sigma_c$ are shown in Figure 15, with an expanded portion on the right. It is apparent that there is no systematic change in the behavior
of the solutions as $\sigma$ increases towards the critical value.
     \begin{figure}[!b]
\begin{center}
\centerline\leavevmode\epsfig{file=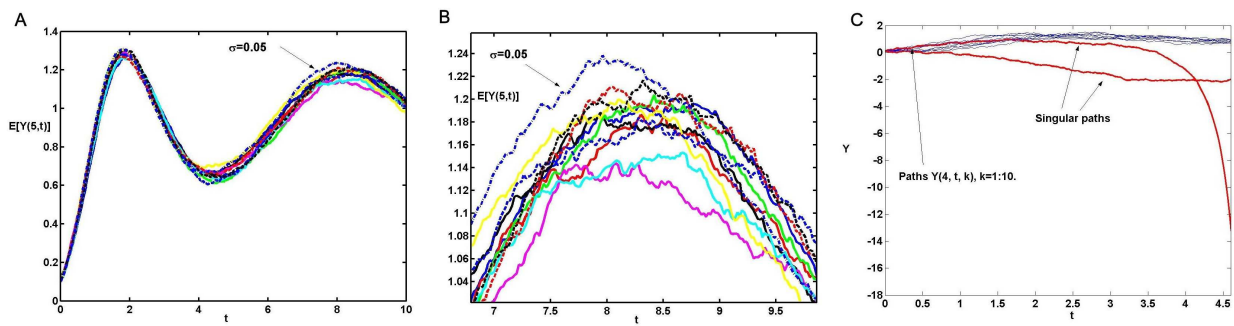,width=6.5in}
\end{center}
\caption{A. Plots of expectations $E [Y(5,t)]$ from simulations  of $Y_{xt}=Y(1-Y) +\sigma W_{xt}$ on $[0,10] \times [0,10]$ for various $\sigma$ 
and initial/boundary value of 0.1  according to Table 3. B.  An expanded version over a small subinterval. C. Ten sample paths at $x=4$ for solutions of
$Y_{xt}=Y(1-Y) +\sigma W_{xt}$ on $[0,10] \times [0,10]$ for  $\sigma=0.1$ 
and initial/boundary value of 0.1. Most sample paths stay 
positive in a narrow band but a few paths (in red), called singular paths,
 take on negative 
values and eventually tend to $-\infty$. 
} 
\label{fig:wedge}
\end{figure}
The 50 sample paths for the smallest value of $\sigma$ at which
the mean became infinitely negative (namely 0.06467148) were analyzed. Divergence to large negative values occurred in just two
of the sample paths.

  Further analysis was made for the case
$\sigma=0.1$ which is considerably greater than  the critical value
$\sigma_c$.  Of ten trials, two resulted in large negative values of $Y$. Trials were stopped when the values of $Y$ became
less than -20 and these are shown in Figure 16, for $x=4$ as a function of $t$
(up to $ t \approx 4.5$ ) along with the 8 trials which 
did not become large and negative. 
It is not known what causes some paths to persist in a downward
trend, although once $Y<0$ the source function $Y(1-Y)$ is negative.
Apparently there were no cases in which paths exceeded unity. 
This property is reminiscent of the anomalous paths in stochastic
Hodgkin-Huxley systems (Tuckwell, 2007, Figure 3).

%
%
%     \begin{figure}[!h]
%\begin{center}
%\centerline\leavevmode\epsfig{file=FQMEAN.eps,width=3in}\epsfig{file=FQSDEV.eps,width=3in}
%\end{center}
%\caption{Mean (left) and standard deviation (right) 
%obtained by  numerical solution (50 trials)
% of $Y_{xt}=Y(1-Y) +\sigma W_{xt}$ with $\sigma = 0.06$
%and initial/boundary value of 0.1.
%} 
%\label{fig:wedge}
%\end{figure}

%
%     \begin{figure}[!h]
%\begin{center}
%\centerline\leavevmode\epsfig{file=FIMSING.eps,width=4in}
%\end{center}
%\caption{Ten sample paths at $x=4$ for solutions of
%$Y_{xt}=Y(1-Y) +\sigma W_{xt}$ on $[0,10] \times [0,10]$ for  $\sigma=0.1$ 
%and initial/boundary value of 0.1. Most sample paths stay 
%positive in a narrow band but a few paths (in red), called singular paths,
% take on negative 
%values and eventually tend to $-\infty$.
%%} 
%\label{fig:wedge}
%\end{figure}

\subsection{Cubic $F$}
The cubic source function is common in simplified modeling in Neuroscience and 
Cardiology where it is used to endow mathematical models
with a resting state and a threshold. Usually the principal
variable (representing a voltage) is coupled to a recovery variable
which is not included here.
Thus the following SPDE is considered briefly
\be Y_{xt}= kY(1-Y)(Y-Y_1) + \sigma W_{xt} \ee
where $k$, $Y_1$ and $\sigma$ are constants.
The values $k=4$ and $Y_1=0.5$ will be used throughout.
\subsubsection{Deterministic example}
For the deterministic case with $\sigma=0$, the equation (31) was
integrated on $[0,30] \times [0,30]$ with a grid $2000 \times 2000$
and with initial ($t=0$) and boundary ($x=0$) values of 0.1 and 
0.6. For both IC/BC values wave solutions formed and attained
apparently constant amplitude and wavelength as $x,t$ became
larger. The values of solutions are plotted as functions of $t$ at $x=15$ 
in Figure 10A.  For both IC/BC values the plots reveal oscillations of about the same period but different amplitudes. For the smaller IC/BC 
value the maximum and minimum of $Y(x,t)$ are 0.1000 and -0.0544,
and oscillations are about the lower equilibrium point $Y=0$.
With the larger IC/BC value the maximum and minimum are 1.1569 and
0.6000 (the IC/BC) value. The mean value of $Y$ is 0.98 as the
oscillations are roughly centered on the upper critical point at $Y=1$.

\subsubsection{Stochastic example}
For the cubic source as above, the effects of the inclusion of noise with $\sigma=0.05$ with an IC/BC value of 0.6 on $ [0,15] \times [0, 15]$
with grid $n_x=n_t=750$ were investigated with 50 trials.   
For the 50 sample paths the minimum value of $Y$ was 0.5996,
just below the IC/BC value of 0.6, and the maximum value was 1.1411,
substantially above the upper critical value of 1. 

For the smaller noise amplitude of $\sigma=0.025$, the mean and standard deviation of $Y$ are plotted in Figures 10B and
10C 
and show that the wavefronts in $E[Y(x,t)]$ and the standard deviation are quite clear with some diffuseness at large $x$ and $t$. 

For the larger value of $\sigma=0.05$, a sample path is
shown in Figure 11A,  revealing discernible but irregular wavefronts.
Figures 11B and 11C give the mean and standard deviation of $Y$  
for this case. 
 The standard deviation starts to climb dramatically at
large $x$ and $t$ values, becoming so large that it flattens the wavefronts of $E[Y]$.  However, it is not known without much more detailed analysis whether this 
enlargement of the standard deviation is a real chance effect 
or due to a deficiency in the numerical integration scheme.

In order to see if using an IC/BC value below the 2nd critical point of $Y=0.5$ could eventually result in a noise induced transition to large amplitude waves, the SPDE was numerically integrated with 
an IC/BC value of 0.45 on  $ [0,15] \times [0, 15]$ with $n_x=n_t=750$
and $\sigma=0.05$. In 50 trials the maximum value of $Y$ was
0.5111 which indicated the absence of large amplitude waves,
When the noise level was increased to $\sigma=0.06$, in 10 trials the
maximum value of $Y$ was 0.4671, but
examination of paths suggested that a large amplitude response might
eventually occur. 
Hence, again with $\sigma=0.06$ the region of integration was increased to $ [0,25] \times [0, 25]$ and a grid with $n_x=n_t=1000$ was employed. In each of ten trials the maximum value of $Y$  exceeded 1.4, well above the upper critical point, indicating that
large amplitude responses had occurred. The maximum value
of the mean was 0.9257 and its minimum value was -0.1435.
Examination of  the ten sample paths showed that after a  few
small amplitude wave-like segments, a sudden transition seemed to
occur to large amplitude wave-like structures. However,
further analysis is required to ascertain whether these transitions
are due to purely stochastic effects.

     \begin{figure}[!h]
\begin{center}
\centerline\leavevmode\epsfig{file=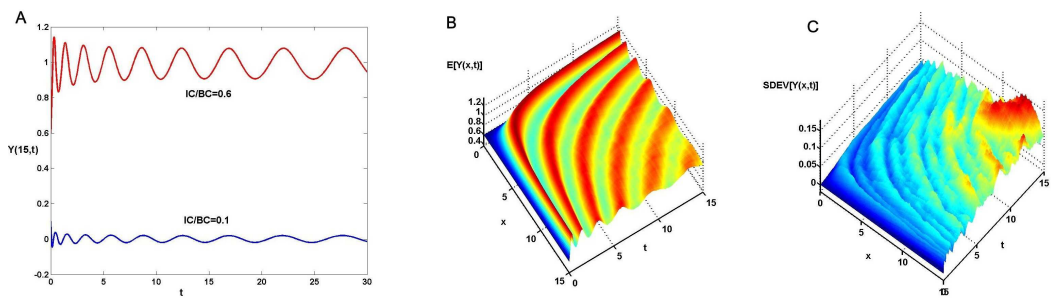,width=5.5in}
\end{center}
\caption{A. Values as functions of time of the numerical solutions of $Y_{xt}=4Y(1-Y)(Y-0.5)$ on
 $ [0,30] \times [0, 30]$ at $x=15$. The blue curve is for the initial/boundary value of 0.1, whereas the red curve is for the initial/boundary value of 0.6.   B and C. Mean and standard deviation of $Y(x,t)$ with the noise parameter
$\sigma=0.025$ and IC/BC value of 0.6.} 
\label{fig:wedge}
\end{figure}

     \begin{figure}[!b]
\begin{center}
\centerline\leavevmode\epsfig{file=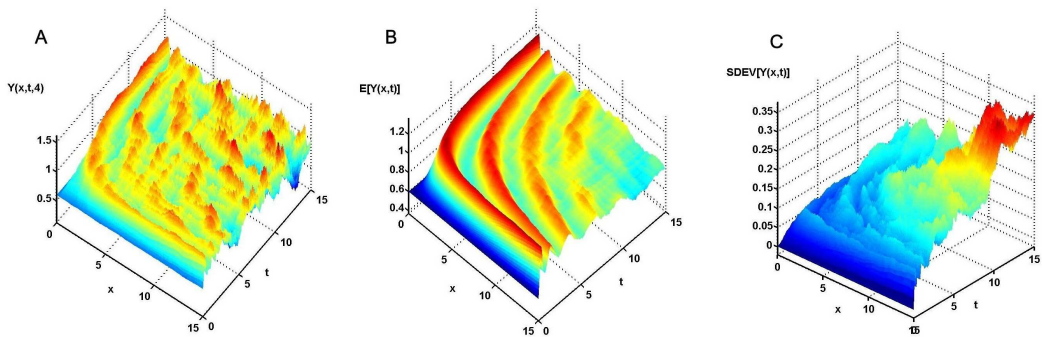,width=5.5in}
\end{center}
\caption{A. A sample path obtained  by numerical integration
of $Y_{xt}=4Y(1-Y)(Y-0.5) + \sigma W_{xt}$ on
 $ [0,15] \times [0, 15]$ with $\sigma=0.05$ and an IC/BC value of
0.6.  B and C. The mean and standard deviation of $Y(x,t)$ from 50 trials.} 
\label{fig:wedge}
\end{figure}

%     \begin{figure}[!t]
%\begin{center}
%\centerline\leavevmode\epsfig{file=FZCUBIC2A.eps,width=3in}\epsfig{file=FZCUBIC2B.eps,width=3in}
%\end{center}
%\caption{The mean and standard deviation of $Y(x,t)$ from 50 trials as described in
%the last figure. 
%} 
%\label{fig:wedge}
%\end{figure}

%     \begin{figure}[!t]
%\begin{center}
%\centerline\leavevmode\epsfig{file=FZCUBIC4A.eps,width=3in}\epsfig{file=FZCUBIC4B.eps,width=3in}
%\end{center}
%\caption{Mean and standard deviation of $Y(x,t)$ for the same
%system as in the previous two figures, but with the noise parameter
%$\sigma=0.025$, being half the prior value. 
%%} 
%\label{fig:wedge}
%\end{figure}

%     \begin{figure}[!t]
%\begin{center}
%\centerline\leavevmode\epsfig{file=FZCUBIC5B.eps,width=3in}\epsfig{file=FZCUBIC5A.eps,width=3in}
%\end{center}
%\caption{Two extended sample paths for 
%$Y_{xt}=4Y(1-Y)(Y-0.5) + \sigma W_{xt}$ on
% $ [0,25] \times [0, 25]$ with $\sigma=0.06$ and an IC/BC value of
%0.45.
%} 
%\label{fig:wedge}
%\end{figure}

\section{Sine-Gordon equation}
The sine-Gordon equation, mentioned above, has a rich and
interesting history cutting across many scientific and mathematical
disciplines.
It evidently made its first appearance in Bour (1862) and Enneper (1870)  in the context of
pseudospherical surfaces and then in Frenkel and Kontorova
 (1939) in the theory of crystal dislocations.

In the 1960s it gained attention as a model for elementary particles starting with the
works of Skyrme (1961) and Perring and Skyrme (1962). Of interest were the soliton solutions called kink and anti-kink. 
By definition, solitons are a special type of solitary wave which maintain their form
after a collision and thus are ascribed particle-like properties.  
The sine-Gordon equation in Euclidean coordinates, 
\be \phi_{tt}-\phi_{xx} = -\sin \phi,  \ee
has the well-known kink/anti-kink (single soliton) solutions
(Marchesoni et al., 1988; Guarcello et al., 2015) 
\be \phi^{K^+, K^-}(x,t) = 4 \arctan \bigg[ \exp \bigg(
  \frac{ \pm (x-x_0-ut) }{\sqrt{1-u^2}} \bigg) \bigg], \ee
where $u$ is speed and $K^{\pm}$ denote kink and anti-kink.
There are also oscillatory soliton solutions called breathers which can be stationary, with frequency $\omega<1$,
\be \phi^{Bs}(x,t) = 4 \arctan \bigg[ \frac{\sqrt{1-\omega^2} \sin \omega t   }{\omega \cosh(x \sqrt{1-\omega^2})}  \bigg], \ee
or moving.

The sine-Gordon equation
 is also satisfied by quantities associated with magnetic fields at Josephson junctions between two superconductors
(Josephson, 1965;  Lebwohl and Stephen, 1967; Scott and Johnson, 1969). More recently the sine-Gordon equation has been mentioned 
in biological settings such as DNA (Yomosa, 1984). Some historical 
review is contained in Braun and Frenkel (1998)  with reference
to the relation between discrete and continuous forms.

% There have been several articles including a number of very recent ones analyzing solutions of the sine-Gordon equation with
%stochastic perturbations as elaborated on below (Pascual and V\'azquez, 1985; Marchesoni, 1986; Marchesoni et al., 1988;
%Biller and Pettrucione, 1990; Hairer and Shen, 2014;
%Anton et al., 2015; Guarcello et al., 2015; Huang et al., 2015). 

\subsection{Some previous studies of stochastic sine-Gordon equations}
One of the first studies of perturbations of sine-Gordon solutions 
was that of Joergensen et al. (1982) in the context of Josephson 
transmission lines.  Included in the
PDE was a term representing a bias current, a loss term and a thermal noise term. The latter was taken to be a 2-parameter Gaussian 
white noise of the type used throughout the present paper.
A similar system was considered for sine-Gordon particles in
Bergman et al. (1983) who put
\be \phi_{tt}-\phi_{xx} = -\sin \phi + \mathcal{F} - G \phi_t \ee
where  $\mathcal{F}$ represents external forces, including noise,
and $G$ is viscosity associated with the loss term. 
Marchesoni (1986) considered the same PDE but without the loss
term in order to see how a random force field might influence kink motion. An assumption that the two-parameter white noise
could be factored into the product of two single-parameter processes
was made and a Langevin equation with a one-parameter noise was analyzed for the relativistic momentum of the soliton solution.
This approach was extended in Marchesoni et al. (1988) to include 
the calculation of nucleation rates. 

Biller and Petruccione (1990) also considered a stochastic sine-Gordon equation of the same form  (35)  as Bergman et al. (1983), but with an additive
noise term which was assumed to have only time-dependence, being a zero-mean 1-parameter Gaussian white noise. They simulated the stochastic PDE and also performed a perturbation analýsis in order to ascertain the effects of random perturbations on the center of mass 
position and speed of the soliton solution.

That approach was not dissimilar to that in another early article
 (Pascual and V\'azquez,1985) on the effects of random perturbations in the sine-Gordon equation. These authors numerically and analytically studied the inclusion of both additive and 
multiplicative white noise. Thus they put 
\be  \phi_{tt}-\phi_{xx} = -\sin \phi + f(\phi, \phi_t; x, t)\ee
 where $f$ is the perturbative term.  Since they were only concerned with small noise
they reduced the problem to seeing how random perturbations
would affect the speed and position of the deterministic soliton solution
(kink or antikink). However, the noise was single-parameter Gaussian white noise which was inserted into a set of ordinary differential equations.  

%In Albeverio et al. (1993) 
Technical analysis of stochastic sine-Gordon equations
and its generalizations has been performed. Recently, such endeavours have been pursued energetically: for example by Hairer and Shen (2015) for 2 space dimensions,  Anton et al. (2015) with two-parameter white noise 
on the unit square, Guarcello et al. (2015) with a numerical
study of breathers in long Josephson junctions and Huang et al. (2015)
with a study of stochastic bifurcation theory.

\subsection{Deterministic example}
The original form of the sine-Gordon equation is 
\be Y_{xt}= \sin Y  \ee
where $x$ and $t$ are sometimes called light-cone coordinates or null coordinates.
 This is in distinction to
the Euclidean coordinates which have been more commonly
used in physics and engineering applications where the more familiar wave equation is as in (25), or more customarily written as 
$Y_{tt}-Y_{xx}=-\sin Y$.  Note that some authors have considered $-\sin Y$ on the right hand side of (37) (for example, Fokas, 1997; Leon and Spire, 2001; Leon, 2003; Pelloni, 2005).

    \begin{figure}[!h]
\begin{center}
\centerline\leavevmode\epsfig{file=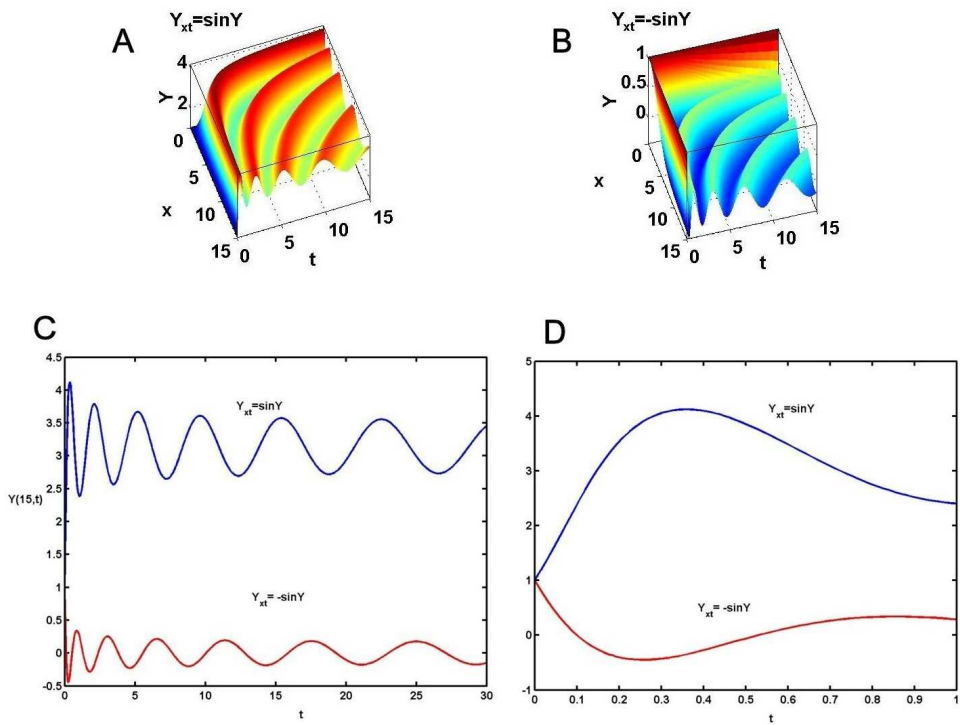,width=5in}
\end{center}
\caption{A. Numerical solution of the deterministic sine-Gordon equation $Y_{xt}=\sin Y$  on  $ [0,15] \times [0, 15]$ with an IC/BC value of
1. B. The corresponding result for 
$Y_{xt}=-\sin Y$. C. Numerically obtained solutions
for $Y_{xt}=\sin Y$  versus time 
at the half-way point $x=15$. D.  
Plot to $t=1$ showing how $Y$ leaves the boundary value
of $Y(15,0)=f(15)=1$.} 
\label{fig:wedge}
\end{figure}

%     \begin{figure}[!h]
%\begin{center}
%\centerline\leavevmode\epsfig{file=FZSINGORD2.eps,width=3in}\epsfig{file=FZSINGORD3.eps,width=3in}
%\end{center}
%\caption{Plots of the numerically obtained solutions in the previous figure versus time 
%at the half-way point $x=15$. In the right hand panel the 
%early parts to $t=1$ show how $Y$ leaves the boundary value
%of $Y(15,0)=f(15)=1$.} 
%\label{fig:wedge}
%\end{figure}
% 

To numerically integrate (37) we use IC/BC conditions as in (12) and (13), and note that existence and uniqueness of solutions for this equation 
with IC/BC of this type on $0 \le x \le \infty$, $0 \le t \le \infty$ were
proved by Angelova et al. (2008) as a Goursat problem. 
The starting point of their analysis was to write the solution as an
integral equation obtained immediately from (37),
\be Y(x,t) =f(x) + g(t)  -c+ \int_0^x\int_0^t \sin Y(u,v)dudv, \ee
where $c=f(0)=g(0)$.
As an example of a deterministic solution we integrate (37)  on $[0,15] \times [0,15]$ with a grid $1000 \times 1000$
and with initial value $f(x)=1$  and boundary value $g(t)=1$.
The computed solution is shown in  Figure 12A. In addition we have
integrated the PDE with $-\sin Y$ on the right hand side using the 
same IC/BC conditions. This solution is shown in Figure 12B. A cross-sectional view of the two solutions at the half-way value of $x=15$  is shown in
Figure 12C. Of interest in the expanded picture of Figure 12D is the smooth way solutions
leave the boundary point at $t=0$.

\subsection{Stochastic examples}
The sine-Gordon equation with additive two-parameter white noise
in the form
\be Y_{xt}=\sin Y + \sigma W_{xt} \ee
was integrated by the numerical method outlined above in Eq. (16).
In all the results to be discussed in this subsection the noise amplitude is $\sigma=0.1$
and the initial/boundary conditions are $f(x)=1$ and $g(t)=1$.

Figures 13A and 13B show the mean $E[Y(x,t)]$ and standard deviation  $SDEV[Y(x,t)]$ based on 50 trials on  $ [0,20] \times [0, 20]$ 
    with a grid 800 x 800.  The mean does not seem to be very 
different from the solution shown in Figure 12A for the noise-free case, although the amplitude of the oscillations in the
mean appears to diminish as the wave progresses.  The standard deviation commences with quite small values at small $x$ and $t$
and begins to steadily grow in magnitude to reach values which are about  4 times larger at (20,20) than (2,2). However, the standard deviation maintains some of the wave-like character of the mean.

     \begin{figure}[!h]
\begin{center}
\centerline\leavevmode\epsfig{file=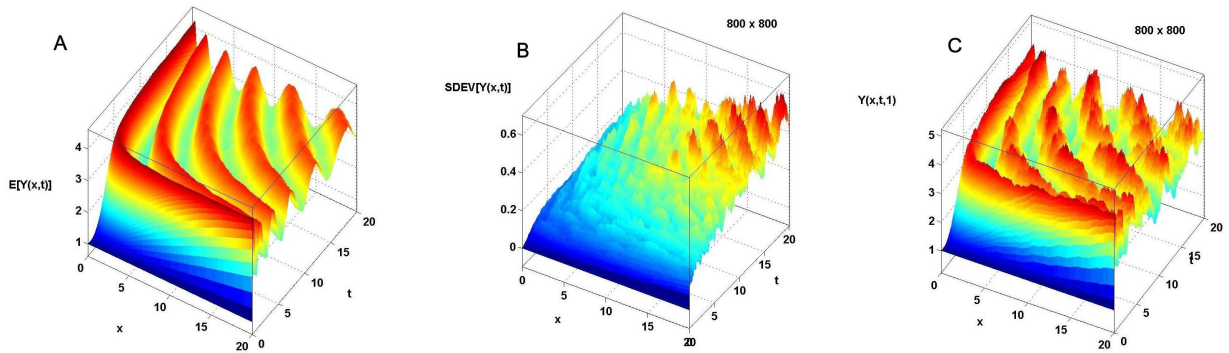,width=6.5in}
\end{center}
\caption{A, B. Mean and standard deviation for numerically
integrated stochastic sine-Gordon equation $Y_{xt}=\sin Y +
\sigma W_{xt}$ with $\sigma=0.1$ on  $ [0,20] \times [0, 20]$ 
    with a grid 800 x 800. Initial/boundary conditions are unity.
C. A sample path for the process as in A,B. } 
\label{fig:wedge}
\end{figure}

In Figure 13C is shown a sample path for $Y$ as in Equ. (32) again
with $\sigma=0.1$ and the same IC/BC as before on  $ [0,20] \times [0, 20]$. 
The wave-like character of the underlying noise-free solution is discernible, which is not expected to be the case with larger values of $\sigma$. 
This was borne out by increasing $\sigma$, firstly to 0.25, 
and then to 0.5 - see Figure 14.

% 
%     \begin{figure}[!h]
%\begin{center}
%\centerline\leavevmode\epsfig{file=SGSAM2800ON10.eps,width=3in}\epsfig{file=SGSAM800ON20.eps,width=3in}
%\end{center}
%\caption{Examples of sample paths for the
%stochastic sine-Gordon equation $Y_{xt}=\sin Y +
%\sigma W_{xt}$ with $\sigma=0.1$
%    and a grid 800 x 800.
%On the left sample path on  $ [0,10] \times [0, 10]$ and 
%on the right another on  $ [0,20] \times [0, 20]$. 
% Initial/boundary conditions are unity.
% } 
%\label{fig:wedge}
%\end{figure}

\subsubsection{Larger values of $\sigma$}

A sample path for $\sigma=0.1$ on $ [0,20] \times [0, 20]$ with 
a grid 1600 x 1600 is shown in Figure 14A. 
For the larger value of $\sigma=0.25$ the sample path  shown
in Figure 14B has a coarse and ragged wave-like structure with
several gaps. The values of $Y$ are mainly positive and do not exceed about 6, which maximum
is slightly larger than that for $\sigma=0.1$.  Thus for $\sigma \le 0.25$
the original wave-like appearance, (on  $ [0,20] \times [0, 20]$)
is roughly preserved.

     \begin{figure}[!h]
\begin{center}
\centerline\leavevmode\epsfig{file=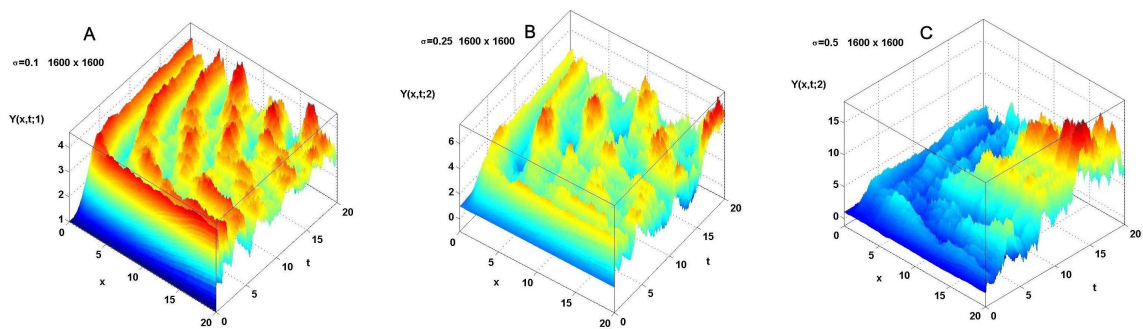,width=6in}
\end{center}
\caption{Sample paths for the
stochastic sine-Gordon equation $Y_{xt}=\sin Y +
\sigma W_{xt}$ 
    and a grid 1600 x 1600 on  $ [0,20] \times [0, 20]$
with $\sigma=0.1$ (A),  $\sigma=0.25$ (B) and $\sigma=0.5$ (C), 
 Initial/boundary conditions are unity.
 } 
\label{fig:wedge}
\end{figure}

With the even larger value of $\sigma=0.5$ the sample
path shown in Figure 14C has a startlingly different character, 
but having a generally wave-like appearance. $Y$   is mainly above zero commencing with a furrowed  L-shaped region (blue) of 
relatively small values which gives way to larger and more haphazard
fluctuations as $(x,t) \rightarrow (20,20)$. 
In another example (not shown),
$Y$ was mainly positive but less than 10, and with a noticeable deep trench which extended down to near -10. In the positive regions there were wave-like structures but
they were not patterned on the deterministic solution 
nor the mean of the small noise case. Thus the wave-like structure was  shattered by the imposed fluctuations.

%
%     \begin{figure}[!h]
%\begin{center}
%\centerline\leavevmode\epsfig{file=SGSIGPT5A.eps,width=3in}\epsfig{file=SGSIGPT5B.eps,width=3in}
%\end{center}
%\caption{Two sample paths for the
%stochastic sine-Gordon equation $Y_{xt}=\sin Y +
%\sigma W_{xt}$ with $\sigma=0.5$
%    and a grid 1600 x 1600 on  $ [0,20] \times [0, 20]$.
% Initial/boundary conditions are unity.
% } 
%\label{fig:wedge}
%\end{figure}

\subsubsection{Increase in number of grid points}
In order to see if making the grid finer has a significant effect on
the numerically calculated solution properties, 50 trials were
executed for the small noise case $\sigma=0.1$  of the sine-Gordon equation (39) with a grid of 
400 x 400 and one of 800 x 800, both on $ [0,10] \times [0, 10]$. 
The standard deviations for these two cases were
not dissimilar, but values with the finer grid were somewhat smaller.  
Wave-like structure was apparent in both results but it seemed to
be better defined with the finer grid. 

%
%     \begin{figure}[!h]
%\begin{center}
%\centerline\leavevmode\epsfig{file=SGSD400ON10.eps,width=3in}\epsfig{file=SGSD800ON10.eps,width=3in}
%\end{center}
%\caption{Comparison of standard deviations of solutions of stochastic
%sine-Gordon equation with $\sigma=0.1$ obtained on  $ [0,10] \times [0, 10]$
%with two different grid sizes: 400 x 400 (left) and 800 x 800 (right).
%Initial/boundary conditions the same at unity. } 
%\label{fig:wedge}
%\end{figure}

%
%
%
\section{Discussion}
As it stands, the form of general equation considered throughout this
article without noise term is
 \be Y_{xt}=F(Y). \ee 
The best-known example is the sine-Gordon equation with $F=\sin Y$ which made its
first appearance over 150 years ago. However, all nonlinear (and linear)
wave equations of the form
\be  \phi_{tt}-\phi_{xx} = G(\phi) \ee can be transformed to the form
of (40) by rotating (and rescaling) the coordinates as for example
in (22) and (23).  

All the numerical solutions in this article were obtained with MATLAB on a PC which puts limits on the mesh/domain sizes through 
limited memory. The question of accuracy of the simulations for  determining properties such as mean and variance for the stochastic cases involves
the fine-ness of the grid, through the magnitudes of $\Delta x$ and $\Delta t$,  and the number of trials. The magnitude of the fluctuations
relative to the solution for $\sigma=0$ is probably also an important factor. The convergence properties of the simple Euler scheme
in the present context are not known, but in the linear case,  when the noise is additive one may compare the mean with the deterministic solution as elaborated on at the start of subsection 3.5. It is likely that for small values of $\sigma$ in the nonlinear examples, the criterion of
approximate agreement of the mean and the deterministic solution
can also be used to gauge the accuracy of the numerical scheme if
the number of trials is large enough. 
Specifically, Figures 6 and 7 for the linear SPDE point to the grid's not being fine enough, as the domain is large,  being  $ [0,30] \times [0, 30]$ 
in Figure 6, and   $ [0,20] \times [0, 20]$ in Figure 7. In Figure
6 the means are compared for various grid sizes but even with the finest grid of 700 x 700 and 50 trials the mean and deterministic solutions differ considerably. Further computations are needed to see whether the disparity diminishes at larger numbers of trials and greater fine-ness of grid.

 For the quadratic SPDE, as seen in Figure 8, agreement of $\hat{Y}$ and 
$Y_{\sigma=0}$ was close, but for the cubic SPDE, Figures 10 and 11
show good agreement only for the first few wavefronts.
For the sine-Gordon equation agreement between $\hat{Y}$ and 
$Y_{\sigma=0}$ is good, as seen in Figure 13 where the domain is  $[0,20] \times [0, 20]$  and the grid 800 x 800 giving $\Delta x = \Delta t =0.025$ which apparently gives good convergence for 50 trials and $\sigma=0.1$.

The boundary conditions employed in the numerical calculations have mainly been constant
$Y$ along $x=0$ and $t=0$, with a consistency requirement at $(0,0)$.
The wave equation (41) is usually endowed with initial conditions
for $u$ and for $u_t$ along with boundary conditions if on a finite space interval. The relation between solutions of (40) and (41) will be explored in a future article. Of much interest will be the determination,
in particular for the sine-Gordon case,  the  nature of solutions in (40) which correspond to soliton solutions in (41). 
 
Another classical PDE of the form of Equ. (40) is Liouville's equation
which can be written as  (Kamran, 2002, who calls $Y_{xt}=F(Y)$ an F-Gordon equation, the case $F=0$ being the wave equation)
\be Y_{xt}=e^Y \ee
with associated wave equation
\be Y_{tt}-Y_{xx}= -e^Y. \ee
  It also has its origins
in the analysis of the Gaussian curvature of a metric in differential geometry.
(Note that there are two other completely different equations called Liouville's; one in statistical mechanics and the other in quantum mechanics). Liouville (1853) gave a general explicit solution of 
equation (42), but it is difficult to consider numerical solutions
of either the deterministic equation or its stochastic forms because
suitable boundary conditions are not available. Preliminary investigations resulted in solutions which rapidly became unbounded. 

Finally, we note that equations of the form of (40)
or more generally
\be Y_{xt}= F(Y, x, t)\ee 
 can constitute a general form of growth model where the space-time rate of change of a biological, chemical or physical quantity depends on population
size as a function of space and time. However, most of the applications of (40) can be, or have been, obtained by the transformation of a wave equation.

% 
%\subsubsection{Summary and conclusions}

\section{Acknowledgements}
%I am grateful to the The University of Adelaide, School of Electrical and Electronic Engineering and its head Professor Chen-Chew Lim for an adjunct position and MIS MPI
%Leipzig through Professor Juergen Jost for an associateship.
This article grew out of preliminary discussions with Professor Fima Klebaner, Department of Mathematics, Monash University.
% to whom I am also grateful for an affiliate position.

\section{References}

\nh Adler, R.J., 1978. 
Some erratic patterns generated
by the planar Wiener process. Suppl. Adv. Appl. Prob. 10, 22-27. 

%\nh Albeverio, S., Di Persio, L.,  Mastrogiacomo, E., 2011.
%Small noise asymptotic expansions for stochastic PDEs. I.  the
%case of a dissipative polynomially bounded non linearity. Tohoku Mathematical
%Journal 63, 877-898.
%
%
%\nh Albeverio, S., Haba, Z., Russo, F. 1993. Stationary solutions of stochastic parabolic and hyperbolic Sine-Gordon equations. Journal of Physics A: Mathematical and General, 26, L711-L718.
%
%\nh Albeverio, S., Mastrogiacomo, E., Smii, B., 2013.
%Small noise asymptotic expansions for stochastic PDEs driven
%by dissipative nonlinearity and L\'evy noise. Stochastic Process. Appl. 123, 2084-2109.  

\nh Angelova, D.T., Georgiev, L.P., Angelov, V.G.  2008.
A Goursat problem for sine-Gordon equation.
Jubilee International Scientific Conference VSU,`2008. 
Available on Google scholar, publisher not stated.

\nh  Anton, R., Cohen, D., Larsson, S.,Wang, X., 2015. Full discretisation of semi-linear stochastic wave equations driven by multiplicative noise. arXiv preprint arXiv:1503.00073.

\nh  Bergman, D. J., Ben-Jacob, E., Imry, Y., Maki, K., 1983. Sine-Gordon solitons: particles obeying relativistic dynamics. Physical Review A, 27, 3345-3348.

\nh Biller. P., Petruccione, F., 1990. Dynamics of sine-Gordon solitons under random perturbations: Multiplicative large-scale white noise. Physical Review B 41, 2139-2144.

\nh Boulakia, M., Genadot, A., Thieullen, M., 2014.
 Simulation of SPDE's for excitable
media using finite elements. hal-01078727.

\nh Bour, E., 1862. Th\'eorie de la d\'eformation des surfaces. Journal
de l'\'Ecole
Polytechnique  39, 99-109.

\nh Braun, O.M., Kivshar, Y.S., 1998. Nonlinear dynamics of the Frenkel-Kontorova model. Physics Reports 306, 1-108.

\nh $\check{\rm C}$encov, N. N., 1956. Wiener random fields depending on several parameters. Dokl. Akad. Nauk SSSR 106 607-09.

\nh Conlon, J.G., Doering, C.R., 2005.
On Travelling Waves for the Stochastic
Fisher-Kolmogorov-Petrovsky-Piscunov Equation.  J. Stat. Phys.
120, 421-477.

\nh D\"orsek, P., Teichmann, J., Velu$\check{\rm s}$$\check{\rm c}$ek, D., 2013. 
Cubature methods for stochastic (partial) differential
equations in weighted spaces. Stoch PDE: Anal. Comp. 1,634-663.

\nh Enneper, A., 1870. \"Uber asymptotische Linien. 
In: Nachr. K\"onigl. Gesellsch. d. Wiss. und G.A. Univ. G\"ottingen 1, 493-510.

\nh Faugeras, O., MacLaurin, J., 2014. Large deviations of an ergodic synchronous
neural network with learning. arXiv:1404.0732v3 [math.PR]. 

\nh Fisher, R.A., 1937.   The wave of advance of advantageous genes.  Ann. Eugen. 7, 355-369.

\nh Fokas, A.S., 1997. A unified transform method for solving linear and certain nonlinear PDEs.  Proc.
R. Soc. A 453, 1411-1443.

\nh Frenkel, J., Kontorova, T., 1939. On the theory of plastic deformation and twinning. Izvestiya Akademii Nauk SSSR, Seriya Fizicheskaya 1, 137-149.

\nh Guarcello, C., Fedorov, K., Valenti, D., Spagnolo, B., Ustinov, A. 2015. Sine-Gordon breathers generation in driven long Josephson junctions. arXiv preprint arXiv:1501.04037.

\nh Hairer, M., Shen, H., 2014. The dynamical sine-Gordon model. arXiv preprint arXiv:1409.5724.

\nh  Hajek, B., 1982. Stochastic equations of hyperbolic type and a two-parameter
Stratonovich calculus. Ann. Probab. 10, 451-463. 

\nh Huang, Q., Xue, C.,Tang, J. (2015). Stochastic D-bifurcation for a damped sine-Gordon equation with noise. AIP Advances 5, 047121.

\nh Joergensen, E., Koshelets, V. P., Monaco, R., Mygind, J., Samuelsen, M. R., Salerno, M., 1982. Thermal fluctuations in resonant motion of fluxons on a Josephson transmission line: Theory and experiment. Physical Review Letters, 49, 1093-1096.

\nh Josephson, B.D., 1965. Supercurrents through barriers.
Adv. Phys. 14, 419-451. 

\nh Kamran, N., 2002. Selected topics in the geometrical study of differential equations (Vol. 96). American Mathematical Society.

\nh Kitagawa, T., 1951.  Analysis of variance applied to function spaces.
 Mem. Fac. Sci. Kyushu Univ. Ser. A. 6, 41-53. 

\nh Khoshnevisan, D., 2001. Five lectures on Brownian Sheet.
http://www.math.utah.edu/˜davar

\nh Khoshnevisan, D., Kim, K., 2015. Nonlinear noise excitation of intermittent
stochastic PDEs and the topology of LCA
groups. Ann. Prob. 43, 1944-1991.  

\nh Kolmogorov, A., Petrovsky, I., Piscunov, N., 1937. 
Study of the diffusion equation with growth of the quantity of matter and its applications to a biological problem. 
Moscou Universitet Bull. Math. 1, 1-25.

\nh Krauss, L.M., 2012.  A Universe from Nothing.  Simon and Schuster, New York.  

\nh Lebwohl, P., Stephen, M.J., 1967. 
Properties of vortex lines in superconducting barriers. Phys. Rev. 163,
376-379.

\nh Leon, J., 2003.  Solution of the Dirichlet boundary value problem for the sine-Gordon equation, Phys. Lett. A 319, 130-142.

\nh Leon, J., Spire, A., 2001. The Zakharov-Shabat spectral problem on the semi-line: Hilbert formulation and applications. Journal of Physics A: Mathematical and General 34, 7359-7380.

\nh Liouville, J. (1853). Sur l'\'equation aux diff\'erences partielles
 ${d^ 2\log \lambda \over du dv} \pm {\lambda \over 2a^ 2}= 0$. Journal de Math\'ematiques Pures et Appliqu\'ees 18, 71-72.

\nh Marchesoni, F., 1986. Solitons in a random field of force: a Langevin equation approach. Physics Letters A 115, 29-32.

\nh Marchesoni, F., H\"anggi, P., Sodano, P., 1988. A Langevin equation approach to sine-Gordon soliton diffusion with application to nucleation rates. In, Universalities in Condensed Matter, pp. 88-92. Springer Berlin Heidelberg.

\nh Pascual, P. J., V\'azquez, L., 1985. Sine-Gordon solitons under weak stochastic perturbations. Phys. Rev. B, 32, 8305-8311.

\nh Pelloni, B., 2005. The asymptotic behavior of the solution of
boundary value problems for the sine-Gordon
equation on a finite interval. J. Nonlinear Mathematical Physics 12, 518-529.

\nh Perring, J.K., Skyrme, T.H.R., 1962.  A model unified 
field equation.  Nucl. Phys. 31, 550-555.

\nh Petterson, K.H., Lind\'en, H., Tetzlaff, T., Einevoll, G.T.,  2014.
Power laws from linear neuronal cable theory: power
spectral densities of the soma potential, soma
membrane current and single-neuron contribution to
the EEG.
PlOS Comp. Biol 10, 1e1003928. 

\nh Scott, A.C., Johnson, W.J., 1969.
Internal flux motion in large josephson junctions.
Appl. Phys. Lett. 14, 316-318. 

\nh Skyrme, T.H.R., 1961. Particle states of a quantized meson field.
Proc. Roy. Soc. Lond. A 262, 237-245.

\nh Stannat, W., 2013. Stability of travelling waves in
stochastic Nagumo equations. arXiv:1301.6378v2 [math.PR]. 

\nh Tuckwell, H.C:, 2007.  Spike trains in a stochastic Hodgkin-Huxley system. BioSystems 80,  25-36.

\nh Tuckwell, H.C., 2008. Analytical and simulation results for the stochastic spatial
Fitzhugh-Nagumo model neuron. Neural Computation 20, 3003-3033.

\nh Tuckwell, H.C., 2013a.  Stochastic partial differential equation models in Neurobiology: linear and nonlinear models for spiking neurons. Springer Lecture Notes in Mathematics 2058,  Stochastic Biomathematical Models, 149-173. Eds. Bachar, M., Batzel, J.,
Ditlevsen, S.  Springer, Berlin. 

\nh Tuckwell, H.C., 2013b.  Stochastic modeling of spreading cortical
depression. Springer Lecture Notes in Mathematics 2058,  Stochastic Biomathematical Models, 187-200. Eds. Bachar, M., Batzel, J.,
Ditlevsen, S.  Springer, Berlin.

\nh Vilenkin, A., 1982. Creation of universes from nothing.
Phys. Lett. 117B, 25-28.

\nh Vilenkin, A., 1984. Quantum creation of universes. Phys. Rev, D. 30,
509-511.

\nh Wellner, J.A., 1975. Monte Carlo of two-dimensional Brownian
sheets. Stat. Inf. Rel. Topics Vol 2, Proc. Summer Res. Inst. Statist.
Inf. Stoch. Proc.  Academic Press, New York. 

\nh Yeh, J, (1960). Wiener measure in a space of functions of two variables, Trans. Amer. Math. Soc. 95, 433-450.

\nh Yomosa, S., 1984. Solitary excitations in deoxyribonucleic acid (DNA) double helices. Phys. Rev. A 30, 474-480. 
%
%Yomosa S 1983 Soliton excitations in deoxyribonucleic acid
%(DNA) double helices; Phys. Rev. A27 2120–2125

\nh Zimmerman, G.J., 1972. 
Some sample function properties of the two-parameter Gaussian process.  
 Annals of Mathematical Statistics 43, 1235-1246.

\end{document}